\documentstyle[12pt]{article}
\begin{document} 
\begin{center}
{\large\bf A DIAGONAL REPRESENTATION OF QUANTUM DENSITY MATRIX USING 
q-BOSON OSCILLATOR COHERENT STATES } 

\vspace{0.5cm}    

R.Parthasarathy and R.Sridhar{\footnote{e-mail addresses: sarathy@imsc.ernet.in ; sridhar@imsc.ernet.in
}} \\

The Institute of Mathematical Sciences \\

C.I.T. Campus, Tharamani, Chennai- 600113, INDIA \\  

\end{center} 

\vspace{1.0cm}

{\noindent{\it{Abstract}}}

\vspace{0.5cm}
 
A q-analogue of Sudarshan's diagonal representation of the Quantum Mechanical 
density matrix is obtained using q-boson coherent states. Earlier result of Mehta and Sudarshan on the
self reproducing property of $\rho(z',z)$ is also generalized and a self-consistent self-reproducing
kernel ${\tilde{K}}(z',z)$ is constructed.

\newpage  

A diagonal representation of the quantum mechanical density matrix using
standard bosonic oscillator coherent states has been given by Sudarshan [1].
 One of the remarkable features of this representation is that for any 
normal ordered operator, its "average expectation value", becomes the same 
as that of the classical function for a probability distribution over the 
complex plane. So, the classical complex representation and the quantum 
mechanical density matrices are in one-to-one correspondence.  

\vspace{0.5cm}
 
The emergence of quantum groups in the study of inverse scattering [2]
resulted in q-Boson [3] and q-Fermion [4] oscillators and many usual quantum 
mechanical descriptions find their appropriate q-analogues [5]. It is 
possible that q-oscillators represent an effective way of dealing with 
non-ideal systems in which interactions are incorporated [6]. It is not 
obvious at present, whether nature makes use of q-Bosons (or q-Fermions) as 
some kind of non-linear excitations of the electromagnetic
(or electron) fields. Nevertheless, it is worthwhile to examine the 
q-anologue of important quantum mechanical descriptions. It is the purpose 
here to find a q-anologue of the diagonal representation of the Quantum 
Mechanical density matrix. Although Nelson and Fields [7] have outlined a 
q-anologue of the diagonal representation of the density matrix, our 
results, while in agreement with them, offer further insight into the 
problem, such as the relationship between q-coherent state matrix elements
and Fock space matrix elements, self-reproducing property and self 
reproducing kernel.

\vspace{0.5cm}

The q-Boson creation  $a^{\dagger}$ and annihilation $a$ operators satisfy 
\begin{eqnarray}
 a a^{\dagger} -qa^{\dagger} a = 1,&;& {[N,a]}  =  -a,\ \  {[N,a^{\dagger}]} = a^{\dagger} ,    
\end{eqnarray}
where $ N(\neq a^{\dagger} a)$ is the number operator. The  q-Boson vacuum state 
$\mid 0>$ is an element of Fock Space ${\cal{F}}$  and 
 $a\mid 0> =0; \  < 0\mid 0> =1$. The $n$ q-Boson normalized 
state $\mid n >$ is given by  
\begin{eqnarray}
\mid n >  =  \frac{1}{\sqrt{[n]!}}(a^{\dagger})^n \mid 0>,  
\end{eqnarray}
where $\{ |n> | n\in {\bf{N}} \cup {0}\}$ and the q-number $ [n] $ is 
\begin {eqnarray}
[n]  =  \frac{1-q^n}{1-q}\  ;\  0\leq q\leq  1 ,
\end{eqnarray}
such that $[n]\rightarrow  n$ as $q\rightarrow   1$. The actions of $a, a^{\dagger}$ and $N$ 
on $\mid n >$  
are  
\begin{eqnarray}
a\mid  n> & =&  \sqrt{[n]}  \mid {n-1}>, \nonumber \\  
a^{\dagger} \mid n > & = & \sqrt{[n+1]} \mid  n+1> ,\ \  N\mid n> =n \mid n>. 
\end{eqnarray}
and $ [n]! = [n][n-1][n-2]....[1].$
The q-Boson coherent state [3,8] is given by 
\begin{eqnarray}
\mid z> & = & \left(e_q^{|z|^2}\right)^{-\frac{1}{2}} 
                \sum_{n=0}^\infty\frac{z^n}{\sqrt{[n]!}}\mid n>,  
\end{eqnarray}
where 
\begin{eqnarray}                 
e_q^x & = & \sum_{n=0}^\infty\frac{x^n}{[n]!}.
\end{eqnarray} 

\vspace{0.5cm} 

In the above $z\in  C$ and $e_q^x$ is the q-exponential function. Using (4) 
it readily follows that $\mid z>$ is an eigen state of $a$ with $z$ 
as eigenvalue, i.e., $ a\mid z> =z\mid z>$. If $\mid z'>$  is another 
such state, then 
\begin{eqnarray}
< z\mid z'>   & = & {( e_q^{|z|^2}\ e_q^{|z'|^2})}^{-\frac{1}{2}} 
 e_q({\bar{z}}z')
\end{eqnarray}
showing such states are {\it{not}} orthogonal. Nevertheless, $< z\mid z> 
=1$. Since $z\in  C$ is continuous and $<z\mid z'>$  
is never vanishing, the states (5) satisfy the criteria of 
coherent states [9]. Furthermore, the states (5) provide a resolution of
unity, namely, [8]
\begin{eqnarray}
\frac{1}{\pi } \int d^2 z \mid z >< z \mid & = & 1.
\end{eqnarray}

\vspace{0.5cm}

Writing $z=r e^{i\theta}$, $  0\leq \theta \leq 2\pi$ , $\theta$ integration is ordinary 
Riemannian type while the $r$-integration will be " q-integration" [10], it 
follows from (5) that 
\begin{eqnarray}
e_q^{r^2}\mid re^{i\theta} ><re^{i\theta}\mid = \sum_{n,m=0}^{\infty}r^{m+n}e^{i(n-m)\theta}
  \frac{1}{\sqrt{[n]![m]!}}\ \mid n><m\mid.  
\end{eqnarray}
Multiplying this by $e^{i\ell \theta}$ and integrating over $\theta$, we find 
\begin{eqnarray}
\int_0^{2\pi}\frac{d\theta}{2\pi}\ e_q^{r^2}e^{i\ell \theta}
\mid r e^{i \theta } < r e^{i \theta}\mid &=&  
    \sum_{n,m=0}^{\infty}\frac{r^{m+n}}{\sqrt{[n]![m]!}}
 \nonumber \\     
             & &          \mid n>< m\mid \delta (\ell +n-m).
\end{eqnarray}

Operating on (10) by the q-differential operator $( \frac {d} {dr} )_q$   
$p$-times, we find
\begin{eqnarray}
{\Bigl( \frac {d} {dr}\Bigr)}_q ^p \int_0^{2 
\pi}\frac{d\theta}{2\pi}\ e_q^{r^2}e^{i\ell \theta}\}\mid  r e^{i \theta }>  
< r e^{i \theta}\mid 
& =  &    \sum_{n,m=0}^{\infty}\frac{r^{m+n-p}}{\sqrt{[n]![m]!}}
 \frac{[m+n]!} {[m+n-p]!}\mid n>< m\mid \nonumber \\  
 & & \delta(\ell +n-m). 
\end{eqnarray}

\vspace{0.5cm}

Now, if we evaluate this equation at $ r=0 $, it is easily seen that 
in the right hand side, the term with $ m+n-p=0 $ alone 
survives. Thus
\begin{eqnarray}
\Bigl\{ {\Bigl( \frac {d}{dr}\Bigr)_q}^p \int_0^{2 \pi}\frac{d\theta}{2\pi}
e_q^{r^2}e^{i\ell \theta}\mid r e^{i \theta }> 
< r e^{i \theta}\mid \Bigr\}_{r=0}   
&=  & \sum_{n,m=0}^{\infty}[m+n]! \nonumber \\  
\frac{1}{[m+n-p]!\sqrt{[n]![m]!}} \mid n>< m\mid \delta(\ell +n-m) 
\delta(n+m-p), \nonumber  
\end{eqnarray}
leading to the result
\begin{eqnarray}
\mid n >< m\mid &  = & 
 \frac{ \sqrt {[n]! [m]!}} {[n+m]!}{\Bigl( \frac {d}{dr}\Bigr)_q} ^{(n+m)} \nonumber \\ 
 & & \int_0^{2 \pi}\frac{d\theta}{2\pi} e_q^{r^2}e^{i(m-n)\theta}\mid r e^{i 
\theta }><r e^{i \theta}\mid {\mid}_{r=0}.
\end{eqnarray}
This result will be used later. A density matrix, in terms of the states in 
$ {\cal {F}} $, has the fpllowing expansion:
\begin{eqnarray}
\rho  &  = & \sum_{n,m=0}^\infty \rho (n,m) \mid n><m \mid .
\end{eqnarray}
Using (12) in (13),we have,
\begin{eqnarray}
\rho  & = &  \sum_{n,m=0}^\infty \rho (n,m)\frac{ \sqrt {[n]! [m]!}} {[n+m]!}{\Bigl( \frac {d}{dr}\Bigr)_q} ^{(n+m)} \nonumber \\ 
 & & \int_0^{2 \pi}\frac{d\theta}{2\pi} e_q^{r^2}e^{i(m-n)\theta}\mid r e^{i 
\theta }><r e^{i \theta}\mid {\mid}_{r=0}.
\end{eqnarray}
which is the desired q-analogue of Sudarshan's diagonal representation of 
the density matrix.

Expansion of $\rho$ in terms of the q-coherent states is 
\begin{eqnarray}
\rho & = & \int d^2z \phi (z) \mid z> <z\mid
\end{eqnarray}
where $ d^2 z = dr^2 d\theta $, with $ dr^2$ integration as q- integration.
The complex function $\phi (z)$ is such that $ \int d^2 z \phi (z) =1$ due 
to $ Tr(\rho )=1$. The expansion coefficients (13) and (15) are related by
\begin{eqnarray}
\rho (n,m) & = & \int d^2 z\   \phi (z) \frac {z^n {\bar z}^m}{\sqrt { [n]! 
[m]!}} {\Bigl(e_q {|z|^2}\Bigr)^{-1}}
 \end{eqnarray}  
which is the q-analogue of Eqn.(7) of Sudarshan [1]. Performing the $\theta$ 
integration,
\begin{eqnarray}
\rho (n,m) & = & \int_0^{{R_q}^2} \phi (r^2){\Bigl\{ \frac {r^{2n}}{[n]!}
{{e_q^{r^2}}^{-1}}}\Bigr\}{\delta}_{n,m}, 
\end{eqnarray}
consistent with $ Tr\rho = 1 $.It is to be noted that the q-exponential 
function necessitates the upper limit as $ R_{q^2} = \frac {1}{(1-q)}$,the
first zero of the q-exponential function [7]. The expression in the curly
bracket in (17) can be interpretted as the q-boson Poisson 
distribution[17].

\vspace{0.5cm}

As a first application of the q-coherent states, we consider the expansion 
of a bounded operator $ F $ as a normal ordered power series in $a$ and  
$a^{\dagger} $, namely,
\begin{eqnarray}
F & = & \sum_{p,q} C_{p,q} {a^{\dagger}}^p {a}^q, 
\end{eqnarray}
with $C_{p,q}$ as expansion coefficients. Using $a \mid z>=z\mid z>$, we have
\begin{eqnarray}
\frac {<z'\mid F \mid z>}{<z' \mid z>} & = & \sum_{p,q} C_(p,q) {{\bar 
{z}}^p }z^q .
 \end{eqnarray}
On the other hand, using (5),
\begin{eqnarray} 
\frac {<z'\mid F \mid z>}{<z' \mid z>} & = &{e_q^{{\bar{z'}}z}}^{-1} 
\sum_{n,m}\ \frac {<n \mid F \mid m>}{\sqrt {[n]![m]!}} (\bar {z'}^n) z^m 
\end{eqnarray}

From (19) and (20),we have 
\begin{eqnarray}
C_{p,q} & = & \sum_{k=0}^{Min(p,q)} \frac {(-1)^k q^{\frac {k(k-1)1}{2}}}
{[k]! sqrt{[q-k]! [p-k]!}} <p-k \mid F \mid q-k>.
\end{eqnarray}
$C_{p,q}$ in Eq(19) can be interpreted as the expansion coefficient of 
the operator $F$ in the q-coherent state space,while $<n\mid F \mid m>$,
the corresponding coefficient in the Fock space and (21) relates them.

As a second application of the representation of $\rho$  interms of 
q-coherent states, we first consider the matrix element of $\rho$ between 
such states. Using (5),(8) and (13),we find
\begin{eqnarray}   
<z'\mid \rho \mid z> & = & \sum_{n,m} \rho (m,n)\Bigl\{e_q^{|z|^2} e_q 
^{|z'|^2}\Bigr\} ^{-\frac{1}{2}}\frac{{\bar{z'}}^m {z}^n}{\sqrt {[m]![n]!}}
\end{eqnarray}
Introduce a function $\rho (z',z)$ as 
\begin{eqnarray} 
\rho (z',z) & = & \sum_{m,n} \rho (m,n) \frac{{\bar{z'}}^m z^n}{\sqrt 
{[m]! [n]!}}
\end{eqnarray}
so that using (22)
\begin{eqnarray}
\rho (z',z) & = & \sqrt{e_q^{|z|^2} e_q^{|z'|^2}}<z' \mid \rho \mid  z>.
\end{eqnarray}
The over completeness of the q-coherent states allows us to rewrite (24) 
as 
\begin{eqnarray}
\rho (z',z) & = & \frac{1}{\pi }\int d^2\xi <z' \mid \rho \mid \xi><\xi 
\mid z>{\Bigl\{ e_q^{|z|^2} e_q^{|z'|^2}\Bigr\}}^{\frac{1}{2}}.
\end{eqnarray}
Substituting for $<\xi \mid z>$ from (7),  (25) can be expressed as
\begin{eqnarray}
\rho (z',z) & = & \int d^2  \xi K(z,\xi)\rho (z',\xi ), 
\end{eqnarray}
where
\begin{eqnarray}
K(z,\xi ) & = & \frac{1}{\pi}{\{e_q^{|\xi |^2}\}}^{-1} e_q^{\bar{\xi}z}.
\end{eqnarray}
Eq(26) attributes the remarkable "self-reproducing property" for 
 $ \rho (z',z)$, with the kernal $K(z,\xi )$  in (27).  Eqns (26) and (27) 
provide the q-analogue of eqns (3.25) and (3.26) of Mehta and 
Sudarshan[11].

From (27), it follows
\begin{eqnarray}
\int d^2 \xi K(z,\xi ) K(\xi ,z') & = & K(z,z'), 
\end{eqnarray}
exhibiting the matrix structure of the kernals. The theory of reproducing 
kernals has been investigated by Aronszajan[12], according to which the 
hermitian nature of $ \rho $ demands the kernal $K(z,\xi )$ to satisfy
\begin{eqnarray}
{K(z,\xi)}^{*} & = & K(\xi ,z). 
\end{eqnarray}
The above kernal (27), as well as the kernal (3.26) of Mehta and Sudarshan do not 
satisfy (29). This can be redressed by finding a suitable reproducing kernal.
Instead of eq(23) for $ \rho (z',z) $,we set
\begin{eqnarray}
{\tilde {\rho}} (z',z)  & = & \sum_{m,n}\rho (m,n) \frac{{\bar{z'}}^m 
z^n}{\sqrt {[m]![n]!}}{\{e_q^{|z|^2} e_q^{|z'|^2}\}}^{-\frac{1}{2}}, 
\end{eqnarray}
which differs from (23) only by the q-exponentials. As these involve $ 
|z|^2 and |z'|^2 $ ,the analyticity of $ \rho (z',z)$ and ${\tilde{\rho}} 
(z',z)$ are the same. With this definition,it follows that 
\begin{eqnarray} 
{\tilde{\rho}}(z',z) & = & \int d^2 \xi {\tilde{K}}(\xi ,z){\tilde{\rho}}(z',
\xi ),
\end{eqnarray}
with 
\begin{eqnarray}
{\tilde{K}}(\xi ,z) & = & \frac{1}{\pi }{\{e_q^{|\xi |^2} e_q^{|z|^2}\} 
}^{\frac{-1}{2}} e_q^{{\bar{\xi }}z}. 
\end{eqnarray}
This reproducing kernal has the desired reproducing property and satisfies 
eq(29), namely the hermiticity requirement. Furthermore, it is seen that 
\begin{eqnarray}
{\tilde{K}}(z,z) \geq 0, 
\end{eqnarray}
and $ {\tilde {K}}(z ,z)$ is an entire function.  

\vspace{0.5cm}

To summarize, we have constructed an explicit diagonal representation 
of the quantum mechanical density matrix
using q-deformed boson coherent states. By expanding the density matrix in terms of q-boson coherent
states ($\in C^{\infty}$) and then in terms of q-boson Fock space states ($\in {\cal{F}}$), the two
expansion coefficients are shown to be related. The matrix element of $\rho$ between q-boson coherent
states facilitated the introduction of $\rho(z',z)$ (23) which has the remarkable self-reproducing
property (26) with respect to a kernel $K(z',z)$ in (27). This acquires significance due to the
non-orthogonality of the q-boson coherent states (7). This provides a q-generalization of the earlier
result of Mehta and Sudarshan [11]. A slightly modified reproducing kernel is introduced which satisfies
all the required mathematical properties of a self-reproducing kernel for density matrix which is
hermitian.           

\vspace{0.5cm}

{\noindent{\bf{References}}}

\vspace{0.5cm}

\begin{enumerate}
\item Sudarshan, E.C.G.: {\it{Equivalence of semiclassical and quantum mechanical descriptions of
statistical light beams}}, Phys.Rev.Lett. {\bf{10}} (1963) 277-279. 
\item Drinfeld,V.G.: {\it{Hopf algebras and the quantum Yang-Baxter Equation}}, Soviet Math. Dokl {\bf{32}} (1985) 254-258 \\
      Jimbo,M.: {\it{A q-difference analogue of $U(g)$ and the Yang-Baxter equation}},  
 Lett.Math.Phys. {\bf{10}} (1985) 63-69.
\item Macfarlane, A.J.: {\it{On q-analogues of the quantum harmonic oscillator and the quantum group
$SU(2)_q$}},   
J.Phys.A: Math.Gen. {\bf{22}} (1989) 4581-4588. \\
      Biedenharn, L.C.: {\it{The quantum group $SU_q(2)$ and a q-analogue of the boson operators}},  
 J.Phys.A: Math.Gen. {\bf{22}} (1989) L873-L878.
\item Parthasaraty, R and Viswanathan, K.S.: {\it{A q-analogue of the supersymmetric oscillator and its
q-super algebra}},  
 J.Phys.A: Math.Gen. {\bf{24}} (1991) 613 - 617.
\item Finkelstein, R.J.:  {\it{Gauged q-Fields}}, hep-th/9906135; 
                         {\it{Observable properties of q-deformed
                         physical systems}},  
                         Lett.Math.Phys. {\bf{49}} (1999) 105-114; 
                         {\it{Spontaneous Symmetry Breaking of q-Gauge Field Theory}}, hep-th/9908210; 
                         {\it{TCP in q-Lorentz theories}}                        
                         Lett.Math.Phys. {\bf{54}} (2000) 157-163; 
                         {\it{On q-Electroweak}}, hep-th/0110075;
                         {\it{On q-SU(3) Gauge Theory}}, hep-th/0206067. \\

  Chaturvedi,S, Jagannathan,R, Srinivasan,V and Sridhar,R.: {\it{Non-relativistic quantum mechanics in a
non-commutative space}},  
J.Phys.A: Math.Gen. {\bf{26}} (1993) L105-L112.  
\item Esteve,J.G, Tejel,C and Villarroya,B.E.: {\it{$SU_q(2)$ quantum group analysis of rotational
spectra of diatomic molecules}},  
      J.Chem.Phys. {\bf{96}} (1992) 5614 - 5617; \\
      Chang,Z, Guo,H-Y and Yan,H.: {\it{The $q$-deformed oscillator model and the vibrational spectra of
      diatomic molecules}},    
      Phys.Lett. {\bf{A156}} (1991) 192 - 196; \\
      Parthasarathy,R.: {\it{$q$-Fermions}}, IMSc-Preprint-93/23, April, 1993.      
\item Nelson,C.A and Fields,M.H.: {\it{Number and Phase uncertainties of the q-analogue quantized field}}.
 Phys.Rev. {\bf{A51}} (1995) 2410 - 2429.
\item Perelemov,A.M.: {\it{On the completeness of some sub systems of q-deformed coherent states}},  
q-alg/9607006.
\item Klauder,R.J and Skagerstam,B.: Eds {\it{Coherent States: Applications in Physics and Mathematical
Physics}}, World Scientific (Singapore), 1985.
\item Exton,M.: {\it{q-Hypergeometric functions and applications}}, Ellis Harwood Limited.
\item Mehta,C.L and Sudarshan,E.C.G.: {\it{Relation between quantum and semiclassical descriptions of
optical coherence}},   
Phys.Rev. {\bf{B138}} (1965) 274-280.
\item Aronszajn,N.: {\it{Theory of reproducing kernels}}, Trans.Am.Math.Soc. {\bf{68}} (1950) 337-404.
\end{enumerate}   
\end{document}